# QUANTIFIER REASONING AND MULTIPLE GENERALITY IN ARISTOTLE AND ANCIENT LOGIC

Clarence Protin, Centro de Filosofia da Universidade de Lisboa (January 2026)

**Abstract**

Aristotelian logic and its related traditions in antiquity are often held to have been equivalent to monadic predicate logic and as such inadequate to formalize mathematics as well as scientific and philosophical discourse in general. In this paper we argue that on the contrary the logical theories of Aristotle and ancient authors such as Galen and Boethius were in fact quite sufficient to account for the logically complex expressions and reasoning involving multiple generality fundamental to the aforementioned disciplines.

## Introduction

A common view is that whereas the deductive mechanisms and logical expressivity of ancient Greek mathematics corresponds to the full power of modern first or higher order predicate logic, it was nonetheless the case that Aristotle's theory of deduction (identified with the theory of syllogism of the Prior Analytics ) captures only monadic predicate logic and was thus radically inadequate for the logical structure and deductive mechanism of mathematics as well as that of scientific and philosophical discourse in general.

In this paper we challenge the above view and argue that the logical theories of Aristotle and some authors of late antiquity (Galen and Boethius) were in fact sufficient to account for the logically complex expressions and deductive processes of ancient mathematics, science and philosophy.

To this end we show first that Aristotle not only makes use of a version of four basic general rules adequate for reasoning about quantified sentences but does so in a philosophical reflected way, that is, the quantifier rules are part of his theory of logic. We also discuss his use of *reductio ad absurdum* and the rules related to implication and inference. We then show that these rules are deployed by Aristotle for sentences containing multiple quantifiers over relations and furthermore present evidence for such a deployment in Galen and Boethius. For Aristotle we were refer to the edition of the Greek text in [topiques1] and [topiques2] for the Topics, [Ros] for the Prior and Posterior Analytics and [Ros2] for the Physics. In the final section we discuss some possible objections to the thesis presented in this paper.

## Quantified sentences and quantifier rules

Quantified sentences are roughly sentences beginning with certain determiners which are commonly called 'quantifiers'. Here we deal only with determiners equivalent either to 'every'/ 'all' / 'for all' (universally quantified sentences) or to 'some' /'there is'/ 'a' (existentially quantified sentences). We take, as did the ancient logicians (among the medieval logicians there were other approaches[pin]), the determiner as affecting the whole sentence rather than affecting primarily the subject. Thus 'all men are mortal' does not mean that the predicate 'mortal' applies to the subject 'all men' but rather says something about how 'mortal' is predicated of 'man'.

Thus quantified sentences are of two basic types with diverse but equivalent syntactic forms.

Universally quantified sentences include the following common forms:

*Every X is P / Every X is such that P holds of it / All X satisfy P/ For all X P holds of it / For all things if they are a X then P holds of them.*

Existentially quantified sentences include the following forms:

*Some X is P / Some X is such that P holds of it / There is a X such that P holds of it / a X is such that P holds of it.*

where X is some term (such as 'man' in the above example) and P is predicate which may have any syntactic structure. Notice that in the constructions headed by a determiner the rest of the sentence is often a phrase involving pronouns and anaphora. We shall return to this when we discuss Boethius' account of complex conditionals.

Non-quantified sentences will be sentences which do not start with above determiners (or with 'There is'). Note that such sentences can still contain quantifiers such as in the conditional : *if all men are mortal then some men are not mortal*.

The four basic quantifier rules are the following.

(Generalized Barbara)[B]

*'For all X P holds of it'*

*'A is X'*

*Therefore 'P holds of A'*.

(Induction)[I]

Suppose that assuming *'A is a X'* we can deduce that *'P holds of A'* without making use of any additional hypothesis beyond *'A is X'*. Then we can deduce (without the above assumption) that *'For all X P holds of it'*. The assumption *'A is a X'* is usually written *'Let A be an X'*.

(note: here we are evidently not using 'induction' in its usual sense of a hypothetical generalization based on a few concrete instances.)

(Some-Introduction)[Si]

*'P holds of A'*

*'A is X'*

Therefore *'Some X is P'*.

(Some-Elimination)[Se]

*'Some X is P'*

Therefore *'A X which is P is P'* and also *'A X which is P is X'*.

And there is a special case:

(Some-Elimination')[Se']

*'A unique X is P'*

Therefore *'The X which is P is P'* and also *'The X which is P is X'*.

We say that a sentence has *multiple generality* if it contains at least one relation (for instance a verb with subject and direct object) in which two distinct arguments are governed by quantifiers. A basic example is 'every man has a father' which is equivalent to 'for any man there is some man such that that man is their father'. Here the pronouns which are the arguments (pronouns with anaphora) of the relation 'father' are governed by two distinct quantifiers.

The quantifier rules above can be considered a variant of the four basic types of natural deduction rule for quantifiers(see [pra]) and (provided we are equipped with adequate rules for propositional reasoning) as equal in strength to modern first and higher order logic. In fact assuming classical propositional reasoning Generalized Barbara and Induction are sufficient to obtain the full power of modern quantifier logic.

The following propositional rules will also figure in our discussion:

(Reductio ad Absurdum)[RA]

If assuming a sentence S and starting from a sentence T we can deduce tha*t not S* then we have that assuming S we can deduce *not T.*

(Modus Ponens)[MP]

If we have that *'if A then B'* and A then we have that B.

(Deduction Rule)[D]

*If under certain assumptions assuming 'A' we can show that 'B' then under the same assumptions we can conclude that 'if A then B'. This rule is so prevalent in Aristotle's proofs that we dispense we further argumentation for its presence.*

Our goal in this paper is to argue that all of the rules above were consciously used by Aristotle and in particular used for sentences with multiple generality (with important evidence provided by Galen).

**The Generalized Barbara rule**

The theory of the syllogism in the Prior Analytics has been the subject of intense investigation in recent decades, specially the modal syllogism. But there is an important overlooked point. A syllogism is supposed to be a rule of deduction and the collection of the syllogisms considered valid is supposed to be a system of rules which can be deployed to form complex proofs or deductions. But what Aristotle does in the Prior Analytics is engage in a series of proofs aimed at showing that certain combinations of sentences are valid syllogisms while others are not. But we can ask: what were the logical rules that Aristotle used in order to prove which syllogisms were valid and which were not ? And to express the syllogisms themselves or things about syllogisms what was the necessary logical or syntactic structure of the sentences employed ?

How is the Barbara syllogism expressed ? *'For all terms A B and C if B is said of all A and C is said of all B then C is said of all A'*. This is patently a universally quantified sentences which involves multiple generality for the relation 'being said of'. And furthermore it is clear that this syllogism considered as a sentence is meant to be used in conjunction with successive applications of the Generalized Barbara rule in order to obtain particular valid instances of the syllogism:

'For all terms A B and C if B is said of all A and C is said of all B then C is said of all A'
'Man is a term'
------------------------------------------------------------------------------------------------------
'For all terms B and C if B is said of all men and C is said of all B then C is said if all men'
'Biped is a term'
------------------------------------------------------------------------------------------------------
'For all term C if biped is said of all men and C is said of all bipeds then C is said of all men'
'Animal is a term'
------------------------------------------------------------------------------------------------------
'if biped is said of all men and animal is said of all bipeds then animal is said of all men'

Consider the discussion of the axioms that need to be applied in a geometric proof in the Prior Analytics 41b13-22. One axiom is

*'for all quantities A B C and D if A = B and C = D then C subtracted from A = D subtracted from B'*

This leaves now doubt about Aristotle's consciousness of universally quantified sentences involving multiple generality (in this case the relation is that of equality) nor about the legitimacy of the Generalized Barbara rule for obtaining instances of the rule.

In fact Galen's *relational syllogisms* [Barn] involve precisely this. For instance Galen [Gal. p.46] mentions:

*'for all quantities A B C if A is the double of B and B the double of C then A is the quadruple of C'*

and Galen's relational syllogism is precisely the Generalized Barbara rule by which we can derive

*'if 8 is double of 4 and 4 is double of 2 then 8 is the quadruple of 2'*

for which we can then further apply Modus Ponens, a rule which was indubitably part of Aristotle's logical theory: in Prior Analytics 53b12-15 we read: *ei gar tou A ontos anagkê to B einai (...) ei oun alêthes esti to A, anagkê to B alêthes einai: if is it necessary that if A is true then B is true and if A is true then necessarily B is true.*

Furthermore in Aristotle's *Topics* we can interpret a *topic* as being a certain kind of universally quantified sentence [prim, slo] (almost always involving multiple generality) which was then instantiated via Generalized Barbara according to the dialectical need at hand, an interpretation which is in accordance to the tradition expounded in Boethius' *De topicis differentiis:* a topic is precisely a 'maximal proposition'. In the Topics a fundamental relation is the genus-species relation and an important topic enunciated in 121b can be written out

*'for all terms A B and C if B is the genus of A and C is the genus of A then either B is contained in C or C is contained in A'*

which again contains multiple generality and is clearly to be instantiated according to the needs of debate using the Generalized Barbara rule.

It is important to note that the multiple generality present in the (possibly multiply) universally quantified topic sentences can be logically quite complex and the quantifications inside the sentence can be for instance over the individuals which belong to the extensions of the universally quantified terms.

**The Induction Rule**

Consider the remarkable passage from chapter 3 of Book II of the Topics (tr. W.A. Pickard):

*Whereas in establishing a statement we ought to secure a preliminary admission that if it belongs in any case whatever, it belongs universally, supposing this claim to be a plausible one. For it is not enough to discuss a single instance in order to show that an attribute belongs universally; e.g. to argue that if the soul of man be immortal, then every soul is immortal, so that a previous admission must be*

*secured that if any soul whatever be immortal, then every soul is immortal. This is not to be done in every case, but only whenever we are not easily able to quote any single argument applying to all cases in common (*euporômen koinon epi pantôn hena logon eipein*) , as (e.g.) the geometrician can argue that the triangle has its angles equal to two right angles.*

What Aristotle is first saying is that in general we do not have 'induction' from the species to the genus. If (essential) accident C applies to species A of B it does not follow that the same accident C applies to B. Thus if the human soul is immortal it does not follow that the soul is immortal (every soul is immortal). Aristotle says that we have to assume that being an accident of the species implies being an accident of the whole genus. Aristotle contrasts this situation to another situation in which *euporômen koinon epi pantôn hena logon eipein* - in which we are able to deduce (in a uniform way) a universal conclusion through a single argument, as the geometrician reasons about property of a triangle. Thus the geometrician may start with 'Let X be a triangle' and arrive at the conclusion that X has its angles equal to two right angles and thereby view this as a proof that all triangles have their angles equal to two right angles. Thus we are lead to interpret Aristotle's *koinon epi pantôn hena logon eipein* as precisely a statement of the Induction rule.

It is also interesting to compare this to Proclus' commentary of the first book of Euclid [proc, p.162] which may be also interpreted as expressing the Induction rule:

*Furthermore, mathematicians are accustomed to draw what is in a way a double conclusion. For when they have shown something to be true of the given figure, they infer that it is true in general, going from the particular to the universal conclusion. Because they do not make use of the particular qualities of the subjects but draw the angle or the straight line in order to place what is given before our eyes, they consider that what they infer about the given angle or straight line can be identically asserted for every similar case. They pass therefore to the universal conclusion in order that we may not suppose that the result is confined to the particular instance. This procedure is justified, since for the demonstration they use the objects set out in the diagram not as these particular figures, but as figures resembling others of the same sort (...)Suppose the given angle is a right angle. If I used its rightness for my demonstration, I should not be able to infer anything about the whole class of rectilinear angles; but if I make no use of its rightness and consider only its rectilinear character, the proposition will apply equally to all angles with rectilinear sides.*

There are a number of proofs in Aristotle's Physics and Organon which make use of variables in the same way as the mathematician's ‘let X be a triangle’ and which cannot be interpreted in any other way than as employing the Induction rule (often in combination with the Deduction rule) and these rules are furthermore applied to sentences containing multiple generality. A typical example is Physics 241b34-242a49, the proof that *hapan to kinoumenon upo tinos anagkê kinesthai – everything that is moved is moved by something.* The beginning and end of the proof read (tr. Hardie and Gaye):

*Everything that is in motion must be moved by something. For if it has not the source of its motion in itself it is evident that it is moved by something other than itself, for there must be something else that moves it. If on the other hand it has the source of its motion in itself, let AB be taken to represent that which is in motion essentially of itself and not in virtue of the fact that something belonging to it is in motion (...) a thing must be moved by something if the fact of something else having ceased from its motion causes it to be at rest. Thus, if this is accepted, everything that is in motion must be moved by something. For AB, which has been taken to represent that which is in motion, must be divisible since everything that is in motion is divisible. Let it be divided, then, at G. Now if GB is not in motion, then AB will not be in motion: for if it is, it is clear that AG would be in motion while BG is at rest, and thus AB cannot be in motion essentially and primarily. But ex hypothesi AB is in motion essentially and primarily. Therefore if GB is not in motion AB will be at rest. But we have agreed that that which is at rest if something else is not in motion must be moved by something. Consequently, everything that is in motion must be moved by something.*

The proof is reduced to the special special case: everything that is changed in its own right is changed by something. Here 'everything' is to be understood as 'every physical body'. The proof involves starting with a variable AB: Let AB be taken to represent, etc. Thus Aristotle begins by employing a variable AB and assuming the hypothesis for AB: 'P holds of AB', where P represents 'changed in its own right'. Then after a detour through a *reductio ad absurdum* argument Aristotle arrives at the conclusion that 'some body moves AB'. After the application of the inference rule to obtain: 'if AP is P then some body moves AB' the conclusion 'every body such that P holds of it is moved by some body' must be seen as following precisely from an application of the Induction rule for the variable AB.

There are passages in the Physics which can be interpreted as involving multiple applications of the Induction rule to derive a conclusion, for instance the theorem in 243a32: the agent of change and that which is changed must be in contact (although there are no variables employed explicitly in the text).

**Some-Introduction**

How does Aristotle show in the *Prior Analytics* that a certain candidate for a syllogism is not valid ? A common strategy involves presenting instances of three terms for which the premises hold but the conclusion patently does not.

It is plausible to assume that Aristotle accepted the rule for converting the negation of universal quantification to existential quantification and the negation of existential quantification to universal quantification for any sentence, not only sentences of a restricted type with only one quantification and without relations. In particular he would accept repeated applications to multiple quantifications.

Thus if he wanted to show that a syllogism was invalid such as in 26a:

*'It is not the case that for all terms A B and C if B is said of all A and C is not said of any B then C is said of all A'*

he would use this conversion rule three times to obtain the equivalent sentence:

*'There are some terms A B and C such that it is not the case that if B is said of all A and C is not said of any B then C is said of all A'*

Note that this sentence again has multiple generality. But an implication is false if the premises are true and the conclusion is false, thus the above sentence is equivalent to:

*'There are some terms A B and C such that B is said of all A, C is not said of any B and C is not said of all A'*

Aristotle gives us the example:

'Man is term', 'Animal is a term', 'Stone is a term', 'Animal is said of all men and stone is not said of any animal and stone is not said of all men'. Then with this example it is implied that we can apply the Some-Introduction rule three times (for 'Man', 'Animal' and 'Stone') to obtain

*'There are some terms A B and C such that B is said of all A, C is not said of any B and C is not said of all A'*

Hence we can assume that Aristotle not only deploys the Some-Introduction rule but does so in the context of multiple generality.

We will present more evidence for the logically conscious use of the Some-Introduction rule in what follows where we discuss the evidence for Some-Elimination, as Aristotle frequently used these two rules together. We note that it would be worth investigating in a future paper the relationship between Some-Introduction and a hypothetical rule which states that from a singular or indefinite sentence (for instance in Boethius' sense which we consider further ahead) we can deduce the corresponding existentially quantified one: from 'It is a man' and 'It is mortal' we deduce 'Some man is mortal'.

If we take 'It is a man' then the deduction 'Something is a man' is taken to employ implicitly the assumption 'It is a thing': we interpret 'Something' as quantifying over a large but nevertheless defined universe of things.

Of interest to the passage between singular or definite sentences to existentially quantified sentences and to the Some-Introduction rule in general is that fact that many of the proofs in Euclid's Elements(cf. [heath]) can be viewed as being based on what is called a *constructivist* interpretation of existential quantification (see for instance [knorr]). The idea is that to prove for instance that 'For all X there is some Y is such that P holds between that X and that Y' amounts to nothing more than producing a construction C which is such that for a given generic A which is X, yields a Y which is in the relation P with X. This can be symbolized as 'Y = C(X)' and 'X and C(X)' are in the P-relation.

Thus in this case the Some-Introduction rule would be read as the following variant:

Given a construction C such that for all A which is X we have that C(X) is a Y which is in a relation P with A then we can deduce that for all X there is some Y which is in the relation P with that X.

**Some-Elimination**

Consider again the proof of the conversion of universal negation in the Prior Analytics 25a14-18 discussed in section 3 (tr. Hardie and Gaye):

*First then take a universal negative with the terms A and B. If no B is A, neither can any A be B. For if some A (say C) were B, it would not be true that no B is A; for C is a B. But if every B is A then some A is B. For if no A were B, then no B could be A. But we assumed that every B is A.*

Aristotle's argument (expressed in a rather terse and laconic style) can be read as follows:

Thesis: If A is not said of any B then neither is B said of any A.

1) Assume A is not said of any B.
2) Assume that B is said of some A (Aristotle takes this immediately to be equivalent to it not being the case that B is not said of any A).
3) Let this some A be c. Then c is a A which is also a B.
4) Hence some B is A.
5) But this contradicts the initial hypothesis that A is not said of any B.
6) Hence by Reductio ad Absurdum we conclude that it is not the case that B is said of some A, that is to say, B is not said of any A.
7) Hence the thesis follows from the Inference Rule applied to 1).

It is important to note that (as generally accepted) Aristotle deploys Reductio ad Absurdim in a logically conscious way: he called this rule *hê eis to adunaton apodeixis* (see for instance Prior Analytics, 41a21-34).

Let us look in more detail at what is going on at 3). We start with ‘some A is B’. Aristotle is saying that c is ‘a A which is B’. But then ‘a A which is B’ is B, that is, ‘c is B’, and also ‘c is A’. But then he deduces immediately that ‘Some B is A’. Thus at 3) it is plausible to interpret the argument in 3) and 4) as involving first the Some-Elimination rule (where ‘a A which is B’ is conveniently renamed c):

From ‘some A is B’ we conclude ‘A A which is B is B’ (i.e. c is B) and ‘A A which is B’ is an A’ (c is A)

followed by the application of Some-Introduction:

From ‘c is A’ and ‘c is B’ (note we switched the order) we conclude ‘Some B is A’.

In 28a24-27 of the Prior Analytics Aristotle gives two proofs of a certain syllogism in the third figure, the second proof employs a process called *tô(i) ekthesthai poiein*, by *ekthesis*. We then confirm that this process corresponds to the deployment of the Some-Elimination rule:

*If they are universal, whenever both P and R belong to S, it follows that P will necessarily belong to some R. For, since the affirmative statement is convertible, S will belong to some R: consequently since P belongs to all S, and S to some R, P must belong to some R: for a syllogism in the first figure is produced. It is possible to demonstrate this also per impossible and by exposition (ekthesis). For if both P and R belong to all S, should one of the Ss, e.g. N, be taken, both P and R will belong to this, and thus P will belong to some R.*

Recall that for Aristotle universal predication has existential import so we can assume that 'some thing is a S'. Aristotle's second proof of this third figure figure (by exposition or ekthesis) can be read in detail as follows:

Assume that 'all S is P' and that 'all S is Q'. We have then that there exists an S (universal predication has existential import, which we read as 'Some thing is S'). Let n be a thing that is S. Then by the hypotheses and using Generalized Barbara we conclude that 'n is Q' and likewise also that 'n is P' . Thus by Some-Introduction 'some P is Q'.

This proof besides presenting strong evidence that proof by ekthesis involved precisely using the Some-Elimination rule (for Alexander of Aphrodisias ekthesis involved an individual which was 'set out' from the extension of some term ) it also presents additional evidence for Aristotle's use of Some-Introduction in the last step: *that the passage from indefinite (or singular) to existentially quantified sentences* was an integral part of his logical theory.

Regarding the 'singular' sentences occurring in the proof it is regrettable that (unlike in Boethius) propositions such as 'Socrates is mortal' are not explicitly discussed in the extant works of Aristotle, specifically with regards to the version of the third-figure: 'Socrates is mortal' and 'Socrates is a man' hence 'Some man is mortal'.

The above proof does not involve multiple generality, but examples which we can interpret as involving the deployment of Some-Elimination to sentences of multiple generality are indeed present in the works of Aristotle.

We will focus on sentences involving multiple generality of the form

*'for all X some Y is such that the X and the Y are in relation R'*

for example 'for every man there is some man that is that man's father'.

If we change ‘there exists’ to ‘there exists a unique’ then we get the modern concept of ‘function’. The use of the Some-Elimination rule for sentences of the above form is found in the Physics.

Typically we have a variable Z (to be used for Induction and introduced as ‘Let Z be such that...’) and we first use Generalized Barbara to instantiate the sentence to

*‘some Y is such that Z and Y are in relation R’*

and then the application of the Some-Elimination rule involves the term *‘a Y such that Z and Y are in relation R’* and yields

*‘Z and a Y such that Z and Y are in relation R are in relation R’* and *‘a Y such that Z and Y are in relation R is a Y’.*

In this case our term *‘a Y such that Z and Y are in relation R’* depends explicitly on Z and it would not be very clear to represent it simply as a constant N; rather is would be better to indicate its dependency by using a functional notation like N(Z). The concept of function corresponds to a form of *genitive*: N(Z) is *Z's Y such that Y is in relation R to it.*

These logical forms abound in the Physics. For example we have ‘every motion has a given time taken’ and therefrom we find ‘let A be the time taken by motion M’. Or more generally ‘for every body and every motion of that body there is a time taken by the motion of that body’ and using this we find expressions of the form ‘let T be the time taken by motion M of body Z’. In the associated proofs in the Physics we thus find the use of Some-Elimination.

One of the questions which present themselves regarding the Topics is how the definition of relations were carried out. While there is evidence that Aristotle did possess such a theory (cf. 145A14-15, 142a26), the details are lacking. To end this section we engage in a speculative reconstruction of how Aristotle would have defined a very basic relation (that of a man being the grandfather of a certain person) and given such a definition how basic reasoning would have been carried out by Aristotle using Some-Elimination (or proof by ekthesis).

There is little doubt that Aristotle would have accepted the definition:

*A being the paternal grandfather of B is A being the father of C and C the father of B for some C*

as well as the axiom

*Every man has a father.* (H)

We will now attempt to reconstruct how Aristotle would have proven that

*Every man has a grandfather.*

A reconstruction of his argument would be something like:

*Let A be a man. Then by (H) A has a father. Let X be A's father. Now again by (H) X has a father. Let B be X's father. Then B is the father of X* (Some-Elimination) *and X is the father of A* (Some-Elimination). *Thus there is some C such that B is the father of C and C is the father of A* (Some-Introduction). *So by definition, B is the paternal grandfather of A. So A has a paternal grandfather. Thus every man has a paternal grandfather (Induction).*

We see here all our familiar rules used together: Induction, Generalized Barbara, Some-Elimination (more rigorously we could use Some-Introduction' as every man has a unique father), Some-Introduction, Deduction and furthermore cases these rules are applied to sentences having multiple generality.

**Multiple Generality in Boethius**

In this section we complement our previous discussion by taking a look a certain passages from Boethius which can be construed as presenting a formal theory of the logical syntax of sentences which we further argue to include multiple generality. We assume that the passages in question represent a transmission of much more ancient material (including that of the Peripatetic schools).

Let us consider passages in Book I 1173D-1176D of Boethius' De topicis differentiis, which are of particular interest. We refer directly to the Latin text [boe].

In 1174D we have a classification of propositions which includes the indefinite and singular type of proposition (the detailed treatment of which is so conspicuously absent in the extant works of Aristotle). The indefinite type of proposition is given by Boethius' example 'Homo iustus est' – 'a man is just'. We will see further ahead the crucial role this type of proposition plays in quantified sentences and multiple generality.

Boethius would accept noun phrases and in particular definite descriptions as subjects. There is evidence for this 1175D:*...partes quas terminos dicimus, non solum in nominibus uerum in orationibus inueniantur* - parts of which are called terms which are found to be not only nouns but also phrases.

An interesting aspect of Book I is Boethius' account of the logical and syntactic structure of propositions. It is clear that Boethius admits what are called 'simple' propositions, involving the standard subject-predicate term relation - and these simple propositions can be universal, particular, indefinite and singular. The terms occurring can themselves be phrases (orationes) - and we already find this in an example of the 'noun' - 'verb' construction given in Plato's Sophist: *Theaetetus, with whom I am now speaking, is flying.* But to Boethius the simple proposition is only one type of proposition. The other kind being the *conditional proposition* which is not simple but complex, being built up from two propositions.

There are two major questions which naturally arise. The first is: is the classification into universal, particular, indefinite and singular valid for propositions in general (thus also for conditional propositions) or only for simple propositions ? The second one is: is the definition of conditional proposition inductive, that is, can the 'parts' of a conditional proposition be any proposition (in particular be themselves a conditional proposition) or are they forced to be simple propositions ?

There is no doubt that the structure of the text – a classification by division - allows us to give an affirmative answer to the first question. After the division into universal, particular, indefinite and singular, Boethius proceeds to divide each of these kinds: 1175A *Harum uero alias praedicatiuas alias conditionales uocamus*. It is of utmost interest to examine the logical and grammatical structure of the examples adduced. But here we confine ourselves to the universal conditional proposition. Note that 'all men are mortal' is expressed ‘for all beings if that being is a man then that being is mortal’ or more succinctly ‘if it is a man it is mortal’. Boethius would say that the antecedent in the simple proposition 'it is a man' and the consequent the simple proposition 'it is mortal'. What kind are such simple propositions? They can only be indefinite. We thus see the key importance of indefinite propositions for the expression of their subjects via personal pronouns with anaphora, *the way they link up subjects from different propositions of a conditional.*

It is Boethius' example of a conditional involving indefinite (or perhaps singular ) simple propositions that allows us to propose an interpretation of sentences expressing multiple generality such as: 'if something happens then something follows it according to a rule'. Boethius would analyze this as having the antecedent ' it happens' and consequent 'something follows it according to a rule'. As for the second question about the definition of conditional being inductive, the answer is indubitably affirmative 1176B: *Harum quoque aliae sunt simplices conditionales aliae coniunctae. Simplices sunt quae praedicatiuas habent propositiones in partibus (...) Coniunctarum uero multiplex differentia, de quibus in his uoluminibus diligentissime perstrinximus, qua de hypotheticis composuimus syllogismis.*

The sentence in 1176A 'Conditionalium uero propositionum, quas Graeci hypotheticas uocant partes, sunt simplices propositiones' is difficult to translate and interpret consistently. Is Boethius saying that simple conditionals are called hypothetical conditionals by the Greeks ? Or that the Greeks call the parts of conditionals hypothetical parts ? But then 'sunt simplices propositiones' does not agree with the discussion on conjunct conditionals which are seemingly also called hypothetical syllogisms.

Only the species of simple conditionals are required to have simple propositions as their parts (their antecedent and consequent) but the same is not true for conjunct conditionals.

Boethius in the passage cited above gives us testimony of the existence of various detailed works on conjunct conditionals : Coniunctarum uero multiplex differentia, de quibus in his uoluminibus diligentissime perstrinximus, qua de hypotheticis composuimus syllogismis.

Where does multiple generality come in ? Recall that the parts of a universal conditional proposition can be any kind of proposition including those with phrase terms. And that it seems that the propositions forming part of the conditional are considered in a special grammatical way via pronouns and anaphora. Consider the sentence: 'if something happens then something follows it according to a rule'. Boethius would analyze this as having the antecedent ' it happens' and consequent 'something follows it according to a rule'. If Boethius accepts in 1176A the phrase predicate ‘investigates the essence of philosophy’ (instantiation of a relation term by a singular term) why would he not, since he admits phrases as predicate or subject terms, accept the predicate ‘being followed by something according to a rule’ ? Further light on this matter might be obtained by analyzing the possible cases of the parts of universal conditional propositions and how the entanglement of pronouns (by anaphora) worked in such cases.

We now present a reconstruction of a fragment of formal logical syntax of the account of propositions given by Boethius in Book I of *De Topicis Differentiis*. This fragment only involves simple conditionals (those involving an antecedent and consequent consisting of simple propositions) and does not include ‘conjunct’ conditionals concerning which, as we saw, Boethius claimed to have written about in detail in several volumes. For simplicity we do not treat negation (or assume that we are considering term-level negation only). The most interesting aspect of this logic is that it includes indefinite and singular propositions, that conditionals can themselves (irrespective of their components) be either universal, particular, indefinite or singular (and a major challenge is to give account of the types of non-universal conditionals) and that entire phrases can function as terms (as subject or predicates).

Let us start with the simple proposition which has the form

S is P

where S and P are ‘terms’. We consider that terms have a ‘valency’ which determines how many arguments they can be applied to (including the subject). Thus in the proposition ‘S is P’ the term P must have valency 1. But the term P (as well as S) can be a ‘phrase’ (oratio). The simplest type of such term (present in the example given by Boethius himself) is given by a verb-term V followed by a direct object term O. The verb-term V must evidently have valency 2 and thus the complex term VO have valency 1 in order to be applied to the subject term S.

But what about the types universal, particular, indefinite and singular of simple propositions. The first two evidently correspond to ‘all S is P’, ‘some S is P’ and the singular proposition ‘S is P’ where S denotes a single object such as Socrates or the sentence ‘1 + 1 = 2’. We can also write singular propositions in the modern form P(S).

Indefinite propositions seem to play a central role in conditional propositions and we argue that they correspond to certain extent to our modern concept of expression with a free variable, the free variable in this case being a third person pronoun (which we assume to be ‘it’). Thus the indefinite simple proposition is ‘it is P’ which we can think of as P(x).

The first hypothesis in our reconstruction involves assuming (and this seems plausible) that Boethius allowed for *abstraction* on propositions which yielded valency 0, 1 and even 2 terms.

That is to say, from a proposition ‘some S is P’ he would allow us to derive the following phrase terms (the cases for ‘all’ and singular propositions are similar):

(being) the case that that ‘some S is P’ (valency 0)
(being) such that some of it is P (valency 1)
(being) such that it holds for some S (valency 1)

For the case ‘some S is VO’ we could further derive the phrase term:

(being) such that some S is V of it (valency 1)

and analogously we could derive from ‘all S is VO’ the phrase term:

(being) such that all S is V of it (valency 1)

Thus Boethius could make use of this last phrase predicate to obtain:

*all S is VO for some O = for some O all S is V of that O*

Consider the example ‘every man knows some things’. Then the proposition ‘some things are known by Socrates’ and then apply the abstraction above to obtain the 1-valent term ‘(being) such that some things are known by it’. We then obtain

*all man is such that some things are known by them*

Also we note that from a simple proposition ‘S is VO’ it is likely that Boethius would allows us to derive the indefinite proposition ‘S is V of it’.

Let us now consider simple conditionals. The form of the universal simple condition is clear. It is of the form if it is A then it is B or, (*for all thing*) if it is A then it is B.

which Boethius considers as having components the indefinite propositions ‘it is A’ and ‘it is B’ which in modern form would be A(x) and B(x) or A(it) and B(it). What is fundamental here is that the two variables are linked grammatically: this corresponds to the anaphora of the pronoun ‘it’ as discussed by Bobzien and Shogry[mult] for this case of Stoic Logic (but for a different set of pronouns corresponding to classical Greek *tis* and *ekeinos*). We could also express our multiple generality example above as the universal simple conditional

if it is a man then some things are known by them

Our challenge now is to give an account of particular, indefinite and singular simple conditionals. We propose for particular conditionals

*something is such that* if it is A then it is B

and for indefinite conditionals (note that these sentences, being themselves indefinite, could be used to further form 'conjunct' conditionals)

*it is such that* if it is A then it is B

With this in mind, we can speculate that an example of a universal conjunct conditional would have been:

(*for any thing*) if it is such that if it is A then it is B then it is C

The most difficult case is that of singular conditionals and indeed it is difficult to see were the standard modern conditional 'if P then Q' comes in.

Our guess for an example of a singular simple conditional would be:

if Socrates is A then Socrates is B

but why not also

if Socrates is A then Plato is B

and indeed the logical form suggests that singular simple conditionals could in general include antecedent and consequents which can either be universal, particular or singular (and thus correspond to standard propositional implication).

**Conclusion**

We have argued for and provided evidence that Aristotle was in possession of a theoretical logic completely adequate for reasoning about quantified sentenced containing multiple generality (together with the basic rules for implication and reductio ad absurdum). This thesis was further corroborated with material from Galen and Boethius. This work can be seen to parallel in the spirit of its historical approach to the work of Bobzien[bob1] and Bobzien and Shogry[mult] on the sophistication of the deductive system of Stoic propositional logic and the adequacy of this logic to express multiple generality. While much remains to be studied and explored regarding quantifier reasoning and multiple generality in medieval logic the conclusion which seems to suggest itself is that our standard view of the history of western logic needs to be somewhat revised so that a greater continuity is acknowledged between ancient logic and the development of quantifier logic by Frege and Peirce as well as the calculus of relations of Schröder. A remarkable consequence is that if our interpretation and reconstruction of Boethius' logic given in this paper be correct then it furnishes some vital clues to better understand Kant's conception of logic as expounded in the Critique of Pure Reason and frees Kant to a certain extent (in so far as he would have been acquainted with Stoic logic and the writings of Boethius) from logical naivety or the imprisonment in first-order monadic logic, a logic unsuitable for mathematics and physics as well as philosophy.

**Objections and Discussion**

Objection: there is a lack of clarity about what qualifies as a logic for Aristotle.

Reply: we take 'logic' as meaning the general theory of deduction, which is precisely the definition of 'syllogism' given in the first book of the Topics.

Objection: just because passages of Aristotle discussing logic contain multiple generality is no reason to conclude that Aristotle's logic encompasses multiple generality.

Reply: we do not adopt this line of reasoning. Rather we have presented evidence from the Analytics, Topics and Physics that Aristotle was aware of the basic natural deduction rules present in mathematical reasoning, rules which are applied to situations of multiple generality in Aristotle's own examples. Furthermore we are charitable to Aristotle in that we deduce from the presence of these rules in proofs and deductions considered as such in the Physics and implicit in the Topics that Aristotle was aware of the formal nature of the deductions present. The paper presents ample evidence for this being the case for both quantifier and propositional rules. We do not argue that syllogistic of the Analytics is capable of handling multiple generality. Rather that this theory is necessarily not to be equated with Aristotle's more general theory of logic (which is present in the deductions of the metatheory of the syllogistic of the Analytics, not in the syllogistic itself). We argue that what can be reconstructed of this logic from the available texts strongly points to Aristotle's general logic being indeed capable of dealing adequately with multiple generality.

Objection: In the Topics Aristotle expresses the various topoi using conditional sentences some of which contain multiple generality. It is a matter of controversy how one uses the topoi in constructing individual syllogisms. But those complex conditional sentences are not themselves deductions: they are assertions about kinds of deductions. There are, moreover, significant obstacles for thinking that the sentences cited in the paper satisfy Aristotle's definition of deduction. Most significantly, the examples cited are complex conditional sentences, but Aristotle seems explicitly to deny that arguments containing conditionals count as deductions.

Reply: we do not argue that topics are themselves deductions. Nor is it correct to classify them as 'complex conditional sentences'. They are before all else universally quantified sentences (as put forward earlier by Slomkovski and Primavesi – and Boethius himself) which need to be instantiated by 'Generalized Barbara' to form (complex) conditionals. The fact that topics are to be used in concrete situations is already a strong argument for the presence of Galen's Relational or Hypothetical

Syllogism in Aristotle.  We also address the question of deductions vs. conditionals in our discussion of Modus Ponens (explicitly enunciated by Aristotle) and Implication Introduction.

Objection:  how do examples from Proclus, Galen and Boethius contribute to your thesis about Aristotle ?

Reply: as the title indicates, this paper takes an organic approach to ancient logic as a whole. Perhaps it would be better to state that we are dealing not only with Aristotle but  with the Peripatetic school as well as well a later syncretistic schools connected to Stoic logic.  We are acutely aware of the scarcity of surviving texts from ancient authors who wrote on logic (just consider the fate of the writings of Chrysippus) and there are many lost writings both of Aristotle and his immediate disciple Theophrastus.  Thus we make the plausible supposition that much of what Boethius, Proclus and Galen wrote could represent the transmission of more ancient material already present in Stoic or Peripatetic schools or even in some lost writings of Aristotle himself. In any case the main thesis of the paper is the adequacy and sophistication of ancient logic in general rather than strictly in Aristotle himself.